\documentclass[11pt]{article}

\usepackage{mathrsfs}
\usepackage{latexsym}
\usepackage{amsmath,amsthm}
\usepackage{amsfonts}
\usepackage{url}

\usepackage{amssymb}
\usepackage{mytheo}

\def\longto{\mathop{\longrightarrow}\limits}

\newcommand{\an}{\alpha_n}
\newcommand{\bn}{\beta_n}

\usepackage{bbm}
\newcommand{\chr}{\boldsymbol{\mathbbm{1}}} 
\newcommand{\pred}[1]{\chr_{\left\{ #1 \right\}}}
\newcommand{\mexp}{\vec{E}}

\newcommand{\E}{\mexp}

\newcommand{\ham}{\textrm{{\tiny \textup{Ham}}}}

\newcommand{\bin}{\operatorname{Bin}}

\usepackage{geometry}
\makeatletter

\newcommand{\ben}{\begin{enumerate}}
\newcommand{\een}{\end{enumerate}}
\newcommand{\bit}{\begin{itemize}}
\newcommand{\eit}{\end{itemize}}

\renewcommand{\vec}[1]{\bs{\mathrm{#1}}}

\newcommand{\dsabs}[1]{\bigl| #1 \bigr|}
\newcommand{\nrm}[1]{\left\Vert #1 \right\Vert}

\newcommand{\R}{\mathbb{R}}
\newcommand{\N}{\mathbb{N}}

\newcommand{\beq}{\begin{eqnarray*}}
\newcommand{\eeq}{\end{eqnarray*}}
\newcommand{\beqn}{\begin{eqnarray}}
\newcommand{\eeqn}{\end{eqnarray}}
\newcommand{\paren}[1]{\left( #1 \right)}

\newcommand{\tlprn}[1]{\left\{ #1 \right\}}
\newcommand{\set}[1]{\tlprn{#1}}
\newcommand{\abs}[1]{\left| #1 \right|}

\newcommand{\floor}[1]{\ensuremath{\left\lfloor#1\right\rfloor}}

\newcommand{\ds}{\displaystyle}

\newcommand{\bs}{\boldsymbol}

\newcommand{\hide}[1]{}
\newcommand{\oo}[1]{\frac{1}{#1}}

\def\eps{\varepsilon}

\newcommand{\pn}{{\hat p}^{(n)}}
\newcommand{\bp}{\boldsymbol{p}}
\newcommand{\bpn}{\hat{\bp}^{(n)}}
\newcommand{\ninf}{\longto_{n\to\infty}}
\newcommand{\bx}{\boldsymbol{x}}
\newcommand{\by}{\boldsymbol{y}}

\title{On 
the
Convergence 
of the Empirical Distribution
}
\author{
Daniel Berend\\
Department of Mathematics and Department of Computer Science\\
Ben-Gurion University\\
Beer Sheva, Israel
\and
Aryeh Kontorovich\\
Department of Computer Science \\
Ben-Gurion University\\
Beer Sheva, Israel}

\begin{document}
\maketitle
\begin{abstract}
We develop a general technique for bounding the tail of the total variation distance
between the empirical and the true distributions over countable sets.
Our methods sharpen a deviation bound of Devroye (1983) for distributions over finite sets,
and also hold for the broader class of distributions with countable support.
We also provide some lower bounds of possible independent interest.
\end{abstract}

\section{Introduction}
Establishing conditions and rates for the convergence of empirical frequencies to their expected values
is a central problem in statistics. 
For concreteness, let $X$ be an $\N$-valued random variable distributed according to $\bp=
(p_1,p_2,\ldots)$
and let $X_1,X_2,\ldots,X_n$ be $n$ independent copies of $X$. The 
canonical estimator for $p_j$ is obtained via the maximum likelihood principle,
which just amounts to a normalized frequency:
\beq
\pn_j = \oo n\sum_{i=1}^n \pred{X_i=j}, \qquad j\in\N.
\eeq
The weak law of large numbers 
guarantees that 
${\ds \pn_j
\ninf
p_j}$ in probability for all $j\in\N$.
The Chernoff-Hoeffding bound $P\paren{\abs{\pn_j-p_j}>\eps}\le2\exp(-2n\eps^2)$, together with the
Borel-Cantelli lemma, strengthens the convergence to be almost sure, thus establishing a strong law of large
numbers. A 
uniform strong law of large numbers
is provided by the Dvoretzky-Kiefer-Wolfowitz inequality
\cite{MR0083864,MR1062069}
\beq
P\paren{\sup_{i\in\N}\abs{\hat F_n(i)-F(i)}>\eps}\le2\exp(-2n\eps^2),\qquad\eps>0,n\in\N,
\eeq
where $\hat F_n(i)=\sum_{j\le i}\pn_j$ and $F(i)=\sum_{j\le i}p_j$. Indeed, since
$\pn_j=\hat F_n(j)-\hat F_n(j-1)$ and $p_j=F(j)-F(j-1)$, we have
\beq
\abs{\pn_j-p_j} &=& \abs{(\hat F_n(j)-\hat F_n(j-1))-(F(j)-F(j-1))}\\
&\le& \abs{\hat F_n(j)-F(j)}+\abs{\hat F_n(j-1)-\hat F_n(j-1)}
\eeq
and therefore
\beq
P\paren{
\nrm{\bpn-\bp}_\infty
>\eps}\le4\exp(-n\eps^2/2),\qquad\eps>0.
\eeq
We conclude that
$
\nrm{
\bpn-\bp
}_\infty 
\ninf
0
$
almost surely (again, Borel-Cantelli is invoked).

An even stronger 
observation
is that $\nrm{\bpn-\bp}_1
\ninf
0$ almost surely.
The $\ell_1$ distance is in some sense 
the most natural one over distributions \cite{gibbs02},
since by Scheff\'e's identity \cite{MR1843146},
\beq
2\sup_{E\subseteq\N}\abs{\bp(E)-\boldsymbol{q}(E)}=\nrm{\bp-\boldsymbol{q}}_1,
\eeq
for any two distributions $\bp,\boldsymbol{q}$ over $\N$ (for
this reason, $\ell_1$ is also referred to as the {\em total variation} distance).
Almost-sure convergence in $\ell_1$
may be surmised from Sanov's theorem \cite{MR2239987,MR1739680} --- whose drawback, however,
is that it does not readily yield explicit, analytically tractable estimates for 
$P\paren{\nrm{\bpn-\bp}_1>\eps}$.

Actually, Sanov's theorem guarantees that
$\bpn\ninf\bp$ in yet a stronger sense, which may be called 
{\em complete convergence in $\ell_1$}. 
Complete convergence was introduced in \cite{MR0019852}.
Applied to the random variable 
\beqn
\label{eq:Jn}
J_n=\nrm{\bpn-\bp}_1,
\eeqn
it means that
\beq
\sum_{n=1}^\infty P(J_n>\eps)
<\infty
\eeq
for all $\eps>0$. 
For $\bp\in\R^k$ 
(that is, distributions with 
support of size $k$),
one may combine the Chernoff-Hoeffding and the union bounds to obtain the following rough
estimate:
\beqn
\label{eq:crude}
P(J_n>\eps) \le 2k\exp(-2n\eps^2/k^2).
\eeqn
Though crude, (\ref{eq:crude}) suffices to 
establish the
complete convergence in $\ell_1$ of $\bpn$ to $\bp$
for distributions with finite support.
A significant improvement is given by \cite[Lemma 3]{MR707939},
which may be stated as follows:
\belen[Devroye]
\label{lem:devroye}
For $\bp\in\R^k$,
we have
\beq
P(J_n
>\eps)\le3\exp(-n\eps^2/25),
\qquad \eps\ge\sqrt{20k/n}.
\eeq
\enlen
However, for 
$\bp\in\R^\N$ with infinite support,
neither
(\ref{eq:crude})
nor
Lemma \ref{lem:devroye}
is applicable.
Our goal in this paper is to establish analogues of Lemma \ref{lem:devroye}
for distributions with countable support. As a by-product, we improve
Devroye's Lemma, sharpening the constant in the exponent by an order of magnitude.

\section{Main results}
\hide{
Our results extend Devroye's Lemma 
in two directions: we give a bound that holds for
all distributions
$\bp\in\R^\N$,
and sharpen the constant in the exponent by an order of magnitude
for distributions with finite support.}
Our basic work-horse is McDiarmid's inequality \cite{mcdiarmid89}, which implies that whenever
$X_i$, $i=1,\ldots,n$, are independent $\N$-valued
random variables and $h:\N\to\R$ is $1$-Lipschitz with respect to the 
Hamming metric\footnote{
The Hamming metric 
is defined by 
$d(\bx,\by) = \sum_{i=1}^n \pred{x_i\neq y_i}$
for
$\bx,\by\in\N^n$.
},
we have
\beqn
\label{eq:mcd}
P(h(X_1,\ldots,X_n)>\E h(X_1,\ldots,X_n)+ n\eps)\le \exp(-2n\eps^2),
\qquad n\in\N,\eps>0.
\eeqn
We choose $h$ to be the function mapping a sample $(X_1,\ldots,X_n)$ to the $\ell_1$ deviation
of 
the empirical frequencies from their expected values:
\beq
h(X_1,\ldots,X_n) &=& \sum_{j\in\N}\abs{ np_j-\sum_{i=1}^n\pred{X_i=j}}
.
\eeq
In the notation above,
$h(X_1,\ldots,X_n)=nJ_n$.
Since $h$ is 
$2$-Lipschitz under the Hamming metric (Lemma \ref{lem:Flip}),
it follows from (\ref{eq:mcd}) that
\beqn
\label{eq:basic}
P(J_n>\E J_n+\eps)\le\exp(-n\eps^2/2),
\qquad n\in\N,~\eps>0.
\eeqn
(In fact, this estimate is near-optimal, as follows from an argument in the spirit of \cite[Theorem 1]{MR1822385}.)

Hence, the crux of the matter is to bound $\E J_n$. For $\bp\in\R^k$, 
it turns out
that
$\E J_n\le\sqrt{k/n}$, which implies our first result:
\bethn
\label{thm:fink}
For every
$k\in\N$,
distribution $\bp\in\R^k$, 
and sample size $n$,
\beq
P(J_n>\eps) \le \exp\paren{-\frac{n}{2}\paren{\eps-\sqrt{\frac{k}{n}}}^2},
\qquad 
\eps\ge \sqrt{\frac{k}{n}}.
\eeq
\enthn
Observe that for 
$\eps\ge\sqrt{20k/n}$, 
Theorem \ref{thm:fink} yields $P(J_n>\eps)\le\exp(-0.3n\eps^2)$, thus 
improving 
Lemma \ref{lem:devroye}.

Our technique works just as well for $\bp\in\R^\N$ with infinite support. Indeed, as we show in
Lemma \ref{lem:nrm1/2},
\beqn
\label{eq:nrm1/2}
\sqrt n \E J_n \le \sum_{j\in\N}\sqrt{p_j} =: \nu(\bp)
,\qquad n\in\N.
\eeqn
When the right-hand side
of (\ref{eq:nrm1/2}) is finite (as is the case for 
``most''
common distributions),
the following result provides a simple and informative bound:
\bethn
\label{thm:infk1/2}
When 
$\nu(\bp)$
is finite, 
\beq
P(J_n>n^{-1/2} \nu(\bp)
+\eps) 
\le \exp(-n\eps^2/2),\qquad n\in\N,~\eps>0.
\eeq
\enthn
When 
$\nu(\bp)$
is infinite, we can still extract meaningful bounds, albeit with a bit more effort.
As we show in Lemma \ref{lem:2sum},
\beqn
\E J_n \le 
\an(\bp)+\bn(\bp)
\label{eq:JAB}
,
\eeqn
where
\beqn
\label{eq:ABn}
\an(\bp) = 2\sum_{p_j<1/n}p_j
,\qquad
\bn(\bp) = 
\oo{\sqrt n}\sum_{p_j\ge1/n} \sqrt{p_j}.
\eeqn
At its most general, our result has the following form:
\bethn
\label{thm:gen}
For all distributions $\bp\in\R^\N$,
\bit
\item[(i)] $P(J_n>\an+\bn+\eps)\le \exp(-n\eps^2/2),\qquad n\in\N,~\eps>0.$
\item[(ii)] $\an+\bn\ninf0$
\item[(iii)] the rate of decay in (ii) may be arbitrarily slow.
\eit
\enthn
The bound in Theorem \ref{thm:gen}(i)
may be rendered effective by our control over $\an$ and $\bn$ for specific distribution families.
Moreover, our estimate 
in (\ref{eq:JAB})
for $\E J_n$ in terms of $\an$ and $\bn$ is nearly tight, in the following sense:
\begin{prop}
\label{prp:ABtight}
For all $n\ge2$ and all distributions $\bp\in\R^\N$, 
\beq
\E J_n\ge \frac{\an+\bn}{4}-\frac{1}{\sqrt n}.
\eeq
\end{prop}
\begin{rem*}
To keep the expressions simple, we have chosen $1/n$ as the break-point
in defining $\an$ and $\bn$.
We note in passing that
a minor improvement in the 
constants 
is achieved by the (optimal)
break-point $1/4n$.
\end{rem*}
The lower bound on $\E J_n$ follows directly from the lemma below, 
in
which 
the first inequality
may be
of independent interest:
\belen
\label{lem:EJbds}
If $Y\sim\bin(n,p)$, then
%
\beq
\sqrt{np(1-p)/2}\le
\E\abs{Y-np}\le \sqrt{np(1-p)},
\qquad
n\ge2,~ p\in[1/n,1-1/n]
.
\eeq
\enlen

\hide{
\section{Notation}
Standard notation is used throughout. In particular, 
$\N=\set{1,2,\ldots}$ denotes the natural numbers and
$\pred{\pi}$ is the 0-1 indicator function
for the predicate $\pi$. We use the vector notation for distributions. 
Thus, $\bp\in\R^k$
or $\bp\in\R^\N$ always consists of nonnegative entries that sum to $1$.
The probability $P(\cdot)$ is always with respect to the product measure
that a specified distribution
$\bp$ induces on $\N^n$.
The $\ell_1$ deviation of the empirical distribution $\bpn$
(obtained from an iid sample of size $n$)
from the true distribution $\bp$ is consistently denoted by
$J_n$ as in (\ref{eq:Jn}).
The {\em Hamming metric} $d_\ham$ on $\N^n$ is defined as follows:
\beqn
\label{eq:dham}
d_\ham(x,y) = \sum_{i=1}^n \pred{x_i\neq y_i},
\qquad x,y\in\N^n.
\eeqn
A function $h:\N\to\R$ is said to be $L$-Lipschitz with respect to $d_\ham$ if
$\abs{h(x)-h(y)}\le Ld_\ham(x,y)$ holds for all $x,y\in\N^n$.
The functions $\an,\bn:\bp\mapsto\R$ are defined as follows:
\beqn
\label{eq:ABn}
\an(\bp) = 2\sum_{p_j<1/n}p_j
,\qquad
\bn(\bp) = 
\oo{\sqrt n}\sum_{p_j\ge1/n} \sqrt{p_j},
\eeqn
and 
\beqn
\label{eq:p1/2def}
\nu(\bp) = \sum_{j\in\N}\sqrt{p_j}.
\eeqn
}
 
\section{Proofs}
We state the following 
elementary
fact without proof:
\belen
\label{lem:Flip}
Suppose $n\in\N$ and $\bp\in\R^\N$ is 
a
distribution.
Define $h:\N^n\to\R$ by
\beq
h(\bx) = \sum_{j\in\N}\abs{np_j-\sum_{i=1}^n\pred{x_i=j}},\qquad \bx\in\N^n.
\eeq
Then $h$ is $2$-Lipschitz with respect to the Hamming metric.
\enlen
\hide{
\bepf
Let the function $\hat n_j:\N^n\to\N$ count the number of times $j$ appears in $x$;
formally, $\hat n_j(x) = \sum_{i=1}^n\pred{x_i=j}$. Now suppose $x,y\in\N^n$
differ only in coordinate $k$, with $x_k=a$ and $y_k=b$.
Then
\beq
{h(x)-h(y)} &=& 
{\sum_{j\in\N}\abs{np_j-\hat n_j(x)} - \sum_{j\in\N}\abs{np_j-\hat n_j(y)}}\\
&=& 
{
\paren{
\abs{np_a-\hat n_a(x)}+\abs{np_b-\hat n_b(x)}
}
-
\paren{
\abs{np_a-\hat n_a(y)}+\abs{np_b-\hat n_b(y)}
}
}\\
&=&
{
\paren{
\abs{np_a-\hat n_a(x)}+\abs{np_b-\hat n_b(x)}
}
-
\paren{
\abs{np_a-(\hat n_a(x)-1)}+\abs{np_b-(\hat n_b(x)+1)}
}
}\\
&\le&
\dsabs{
\abs{np_a-\hat n_a(x)}
-
\abs{np_a-(\hat n_a(x)-1)}
}
+
\dsabs{
\abs{np_b-\hat n_b(x)}
-
\abs{np_b-(\hat n_b(x)+1)}
}\\
&\le&2.
\eeq
\enpf
}

\belen
\label{lem:nrm1/2}
Suppose $n\in\N$ and $\bp\in\R^\N$ is 
a
distribution.
Then
\beq
\sqrt n \E J_n \le 
\sum_{j\in\N}\sqrt{p_j} 
.
\eeq
\enlen
\bepf
Let $Y_j
\sim
\bin(n,p_j)$.
Then
\beq
\paren{\E\abs{Y_j-np_j}}^2 \le \E(Y_j-np_j)^2 
= np_j(1-p_j) \le np_j,
\eeq
whence
\beqn
\label{eq:Ynp}
\E\abs{Y_j-np_j}
\le\sqrt{np_j(1-p_j)}
\le\sqrt{np_j}.
\eeqn
Since
\beqn
\label{eq:nJn}
n\E J_n = \sum_{j\in\N}\E\abs{Y_j-np_j},
\eeqn
the claim follows.
\enpf

\belen
Suppose $n\in\N$ and $\bp\in\R^\N$ is 
a
distribution.
Then
\label{lem:2sum}
\beq
\E J_n 
&\le& \an+\bn.
\eeq
\enlen
\bepf
As in the proof of Lemma \ref{lem:nrm1/2}, 
let
$Y_j\sim\bin(n,p_j)$ and 
use (\ref{eq:nJn}) to obtain
\beqn
n\E J_n 
&=& 
\label{eq:EJAB}
\sum_{p_j<1/n}\E\abs{Y_j-np_j} 
+
\sum_{p_j\ge1/n}\E\abs{Y_j-np_j}.
\eeqn
By (\ref{eq:Ynp}),
the second term on the right-hand side of (\ref{eq:EJAB})
is clearly upper-bounded by $n\bn(\bp)$.
To bound the first term, we appeal to the {\em mean absolute deviation}
formula for the binomial distribution \cite{kenney62}
\beqn
\label{eq:MAD}
\E\abs{Y_j-np_j}=
2(1-p_j)^{n-\floor{np_j}} p_j^{\floor{np_j}+1}\paren{\floor{np_j}+1}\binom{n}{\floor{np_j}+1},
\eeqn
which simplifies to
\beqn
\label{eq:MADappr}
\E\abs{Y_j-np_j}= 2n(1-p_j)^np_j \le 2np_j,
\qquad p_j<1/n.
\eeqn 
This shows that the first term on the right-hand side of (\ref{eq:EJAB})
is upper-bounded by $n\an(\bp)$ and proves the claim.
\enpf

\bepf[Proof of Theorem \ref{thm:fink}]
We claim that
\beqn
\label{eq:EJn}
\E J_n \le \sqrt{\frac{k}{n}}.
\eeqn
Indeed, 
by
Lemma \ref{lem:nrm1/2},
\beqn
\label{eq:sqrtnEjn}
\sqrt n \E J_n \le \sum_{j=1}^k\sqrt{p_j} 
.
\eeqn
Define $\bx\in\R^k$ by $x_j=\sqrt{p_j}$ and recall that
\beqn
\label{eq:nrmx}
\sum_{j=1}^k\sqrt{p_j}=\nrm{x}_1\le\sqrt k\nrm{x}_2=\sqrt k,
\eeqn
which yields (\ref{eq:EJn}).
In view of
(\ref{eq:basic}), this implies the theorem.
\enpf

\bepf[Proof of Theorem \ref{thm:infk1/2}]
Immediate from 
(\ref{eq:basic})
and
Lemma \ref{lem:nrm1/2}.
\enpf

\belen
\label{lem:AB}
Let $\bp\in\R^\N$ be 
a
distribution.
Then
\beq
\an(\bp)+\bn(\bp)\ninf0.
\eeq
\enlen
\bepf
The decay of $\an(\bp)$ to zero is obvious, since
it is the tail of a convergent series. To prove 
that 
\beqn
\label{eq:Bclaim}
\lim_{n\to\infty} \oo{\sqrt n}\sum_{p_j\ge1/n} \sqrt{p_j} =0,
\eeqn
we define the function $\sigma:\N\to2^\N$ by
\beq
\sigma(n) = {\set{j\in\N:p_j\ge1/n}}.
\eeq
Since 
as in (\ref{eq:nrmx}),
\beq
\sum_{p_j\ge1/n} \sqrt{p_j} \le \sqrt{\abs{\sigma(n)}},
\eeq
it suffices to show that
\beq
\abs{\sigma(n)} = o(n).
\eeq
Suppose, to the contrary, that there exist a $c>0$ and an increasing
sequence $(n_k)_{k=1}^\infty$ such that
\beq
\abs{\sigma(n)} \ge
cn_k,\qquad k\ge 1.
\eeq
Put $n_0=1$. Passing to a subsequence, we may assume that
$n_k\ge 2n_{k-1}/c$ for every $k\ge1$. Now
\beq
1&=&\sum_{j=1}^\infty p_j \\
&\ge& \sum_{k=1}^\infty \sum_{\oo{n_{k}}\le p_j <\oo{n_{k-1}}} p_j \\
&\ge& \sum_{k=1}^\infty 
\paren{\abs{\sigma(n_{k})}-\abs{\sigma(n_{k-1})}}
\cdot\oo{n_k}\\
&\ge& \sum_{k=1}^\infty (cn_k-n_{k-1})\cdot\oo{n_k}\\
&\ge& \sum_{k=0}^\infty (cn_k-cn_{k}/2)\cdot\oo{n_k}
= \sum_{k=0}^\infty \frac{c}{2}=\infty.
\eeq
The contradiction completes the proof.
\hide{
Since obviously $\abs{\sigma(n)}\le n$,
the negation of (\ref{eq:sigmaon}) implies
the existence of $0<c<\infty$ such that
\beqn
\label{eq:sigmaCn}
cn \le \abs{\sigma(n)} \le n
\eeqn
holds for infinitely many $n\in 
\N$. 
Let us suppose (to get a contradiction) that (\ref{eq:sigmaon}) fails to hold
and record some basic observations:
\bit
\item[(i)]
for $n_1<n_2$, we have $\sigma(n_1)\subset\sigma(n_2)$
\item[(ii)] $\forall m\in\N~ \exists n_2>n_1>m$ such that
\beq
\abs{\sigma(n_1)} \le n_1 ,\qquad
2n_1 \le cn_2 \le \abs{\sigma(n_2)}
\eeq
\item[(iii)] for any infinite increasing sequence $\set{n_s}\subset\N$, we have
\beq
\sum_{j\in\N} p_j &=& \sum_{s\in\N} 
\nrm{ \sigma(n_{s+1})\setminus\sigma(n_{s})}
,
\eeq
where 
$\nrm{E}$ denotes $\sum_{x\in E}x$ for finite $E\subset\R$.
\eit
Pick any $m\in\N$ and let $n_1,n_2$ be as stipulated in (ii). 
Then
\beq
\nrm{ \sigma(n_2)\setminus\sigma(n_1)} \ge \frac{cn_2}{2}\cdot\oo{n_2}=\frac{c}{2},
\eeq
since the set 
$\sigma(n_2)\setminus\sigma(n_1)$ contains at least $cn_2/2$ elements of magnitude at least $1/n_2$.
Taking $\set{n_s}\subset\N$ to be the infinite sequence 
satisfying (\ref{eq:sigmaCn}),
we have that $\sum_{j} p_j$ diverges by Cauchy's criterion. This contradicts $\sum_jp_j=1$,
and so (\ref{eq:sigmaCn}) can only hold for finitely many $n$;
thus (\ref{eq:sigmaon}) is proved.
}
\enpf

\belen
\label{lem:rates}
For 
any
rate sequence $1>r_1>r_2>\ldots\searrow0$,
there is a distribution $\bp\in\R^\N$ such that 
\beq
\an(\bp)+\bn(\bp)
>
r_n, \qquad n\in\N
.
\eeq
\enlen
\bepf
It suffices to show that there is no rate sequence bounding $\an$.
But this is obvious, since $\an$ may be expressed as the tail of a
series converging to $2$ --- and although any such tail must decay to zero,
the rate may be arbitrarily slow. 
In particular, given some rate sequence $(r_n)$, to ensure
that $\sum_{p_j\ge1/n}p_j\le 1-r_n$ for each $n\in\N$, we may choose the appropriate $p_j$
in an iterative greedy fashion, for $n=1,2,\ldots$.
\enpf

\bepf[Proof of Theorem \ref{thm:gen}]
Item (i) is an immediate consequence of (\ref{eq:basic}) and (\ref{eq:JAB}).
Items (ii) and (iii) are the contents of Lemmas \ref{lem:AB} and \ref{lem:rates}, respectively.
\enpf

\bepf[Proof of Lemma \ref{lem:EJbds}]
The upper bound is contained in (\ref{eq:Ynp}) --- and in fact, holds for all $p$.
To establish the lower bound,
let us rewrite the mean absolute deviation formula
(\ref{eq:MAD}) as
\beq
\E\abs{Y-np}=2k{n\choose k}p^k(1-p)^{n-k+1},
\qquad (k=\floor{np}+1)
.
\eeq
Denote the right-hand side by $E(n,k,p)$, and put
$G(n,k,p)=2E(n,k,p)^2/(p(1-p))$.
The left-hand inequality in the lemma 
is equivalent to the claim
\beqn
\label{eq:Gnkp}
G(n,k,p)\ge n,
\qquad
p\in[1/n,1-1/n],
~
k=\floor{np}+1
.
\eeqn
The domain where (\ref{eq:Gnkp}) is to be proved may be reparametrized
by the inequalities
\beq
2\le k\le n-1,\qquad
\frac{k-1}{n}\le p < \frac{k}{n}.
\eeq
Now
the function
$G(n,k,\cdot)$ is 
increasing on $[(k-1)/n,(2k-1)/2n]$
and decreasing on $[(2k-1)/2n,k/n]$ --- and hence
we need only consider the
endpoints $p=(k-1)/n$ and $p=k/n$.

To examine the first possibility, we take $p=(k-1)/n$ and seek a $k$ that minimizes
$G(n,k,(k-1)/n)$. To 
this end,
we consider the inequality
$G(n,k+1,k/n)\ge G(n,k,(k-1)/n)$,
which is 
equivalent (after a routine calculation)
to
\beqn
\label{eq:nk}
\paren{\frac{k}{k-1}}^{2k-1}\ge
\paren{\frac{n-k+1}{n-k}}^{2n-2k+1}.
\eeqn

Since the function $f(x)=(1+1/x)^{2x+1}$ is monotonically decreasing on $[1,\infty)$,
the inequality
(\ref{eq:nk}) holds whenever $k\le(n+1)/2$. We conclude that
$G(n,k,(k-1)/n)$ 
is minimized
at the smallest allowed value of $k$, which is $k=2$.
We easily verify that the inequality
$G(n,2,1/2)\ge n$
is equivalent to 
$8(n-1)^{2n-1}\ge n^{2n-1}$ for all $n\ge2$, which again follows from the
monotonicity of $(1+1/x)^{2x+1}$. 

The second case, $p=k/n$, is analyzed in an exactly analogous manner.
\enpf

\bepf[Proof of Proposition \ref{prp:ABtight}]
Let $n\ge2$ and
$Y_{j}\sim\bin(n,p_j)$.
We 
group the probabilities 
as follows:
$S_1=\set{j:p_j<1/n}$,
$S_2=\set{j:1/n\le p_j \le 1/2}$
and
$S_3=\set{j:p_j>1/2}$.
By (\ref{eq:MADappr}) and 
Lemma \ref{lem:EJbds},
\beq
\E\abs{Y_{j}-np_j}&\ge&
\oo2
\left\{
\begin{array}{ll}
np_j,& j\in S_1
\\
\sqrt{np_j},& j\in S_2
\end{array}\right.
\hide{
\E\abs{Y_{j}-np_j}&\ge&np_j/2
,\qquad j\in S_1 \\
\E\abs{Y_{j}-np_j}&\ge& \sqrt{np_j}/2
, \qquad j\in S_2
}
.
\eeq
Now
\beq
n\an(\bp) = \sum_{j\in S_1}2np_j 
\le 4\sum_{j\in S_1}\E\abs{Y_{j}-np_j}
\eeq
and
\beq
n\bn(\bp) &=& \sum_{j:p_j\ge1/n}\sqrt{np_j} \\
&\le& \sum_{j\in S_2}\sqrt{np_j} + \sqrt{n} \\
&\le& 2\sum_{j\in S_2}\E\abs{Y_{j}-np_j} + \sqrt{n}
\eeq
and thus
\beq
4\sum_{j\in S_1}\E\abs{Y_{j}-np_j} +
2\sum_{j\in S_2}\E\abs{Y_{j}-np_j} + \sqrt{n}
&\ge& n\an+n\bn,
%
\eeq
which proves the claim.
\hide{
For 
$j\in S_1$,
(\ref{eq:MADappr}) implies
\beq
(2/e)np_j 
\le 
\E\abs{Y_{j}-np_j}\le 2np_j,
\eeq
whereas for 
$j\in S_2$,
Lemma \ref{lem:EJbds} yields
\beq
\oo2\sqrt{np_j}\le\E\abs{Y_{j}-np_j}\le\sqrt{np_j}.
\eeq
}

\hide{
Finally, note that $S_3$ can contain at most one element
and thus $\sum_{j\in S_3} \E\abs{Y_{j}-np_j}\le \sqrt{n}$.
It follows that
\beq
n\E J_n &=& \sum_{j\in\N}\E\abs{Y_{p_j}-np_j} \\
&\ge& \sum_{j\in S_1}(2/e)np_j + \sum_{j\in S_2}
\eeq
}
\hide{
The only 
case 
unaccounted for 
is $p\ge1-1/n$, 
which can contribute
at most $2$ to $n\E J_n$ (note that (\ref{eq:MAD}) is symmetric in $p$
and $1-p$).
}
\enpf

\section*{Acknowledgements}
We thank to Larry Wasserman for referring us to
Devroye's Lemma,
and David McAllester for reminding us about Sanov's Theorem.
We are grateful to the Stone family for providing a venue for this work.

\bibliographystyle{plain}
\bibliography{../mybib}

\end{document}